\documentclass[12pt]{amsart}
\begin{document}

\author{S. Adlaj}
\address{CC RAS, Vavilov st. 40, Moscow, Russia, 119333}
\email{SemjonAdlaj@gmail.com}
\title{A point of order 8}
\begin{abstract}
A formula expressing a point of order 8 on an elliptic curve, in terms of the roots of the associated cubic polynomial, is given. Doubling such a point yields a point of order 4 distinct from the well-known points of order 4 given in standard references such as ``A course of Modern Analysis'' by Whittaker and Watson.
\end{abstract}
\date{May 30, 2011}
\maketitle
\parskip 6pt
\thispagestyle{empty}

\noindent
Let
\begin{equation}
y^2 = 4(x-e_1)(x-e_2)(x-e_3), e_i \in \mathbb{C}, \ i=1,2,3, \label{DefiningEquation}
\end{equation}
be the defining equation for an elliptic curve E over the complex field $\mathbb{C}$, and let $$\beta:= \sqrt{\frac{e_1 - e_3}{e_1 - e_2}}, \ \gamma:= \sqrt{(e_1 - e_3)(e_1 - e_2)}, \ i:= \sqrt{-1}.$$
The roots of the cubic on the right hand side of the defining equation (\ref{DefiningEquation}) need not sum to zero but assume that $\beta>1$, and introduce the values
$$\beta_1:= {\sqrt{\frac{\beta + 1}{\beta - 1}} + \sqrt{\frac{2}{\beta - 1}} > \beta_2:= \sqrt{\frac{2}{\beta + 1}} + \frac{1}{\beta}}.$$
The point $P = (x,y)$ on $E$, where
$$x = e_1 - \gamma - \gamma \left( \sqrt{\frac{\beta + 1}{2}} - 1 \right)$$
$$\left( 1 - \frac{1}{\beta} + \sqrt{1 + \frac{1}{\beta}} \left( \sqrt{\beta_1 + \beta_2} + i \left( \sqrt{\beta_1 - \beta_2} + \sqrt{1 - \frac{1}{\beta}} \ \right) \right) \right),$$
is then a point of order 8.

\noindent
Note that doubly doubling the afore-indicated point $P$ cannot possibly yield either $(e_1,0)$ nor $(e_3,0)$, and so must, if $P$ is indeed of order $8$, yield the point $(e_2,0)$.

\noindent
For an example, consider an elliptic curve E given in Weierstrass normal form via \mbox{equation (\ref{DefiningEquation})}, where
$$e_1 = i, \ e_2 = 0, \ e_3 = -i.$$
Then
$$\gamma = i \sqrt{2}, \ \beta = \sqrt{2}, \ \beta_1 = 1 + \sqrt{2} + \sqrt{2 \left( \sqrt{2} + 1 \right)}, \ \beta_2 = \frac{1}{\sqrt{2}} + \sqrt{2 \left( \sqrt{2} - 1 \right)},$$
and the $x$--coordinate of $P$ is calculated to be
$$x(P) = \sqrt{2} - 1 - i \sqrt{2 \left( \sqrt{2} - 1 \right)}.$$
One might employ the well-known doubling formula, found in standard sources such as \cite{WW}, for successively calculating the points $2P = (1, \pm 2 \sqrt{2})$ and $4P = (0,0)$, the latter evidently being a point of order 2 on E.

\noindent
Incidentally, the formulas given in \cite[\textsection 20.33, p. 444]{WW} apply here yielding two pairs of points of order 4, whose $x$-coordinates are $e_1 \pm \gamma = i \hspace{.04cm} (1 \pm \sqrt{2})$, so all four points differ from either point of the already computed pair $2P$, as they match when doubled the points $(\pm i,0)$, which are, aside from the point $4P$, the two remaining points of order 2 on E.

\end{document}